\documentclass[12pt]{article}
\usepackage[T1]{fontenc}
\usepackage[latin1]{inputenc}
\usepackage{amsmath,amsfonts,amssymb,amscd,mathrsfs,latexsym,epsfig}

\usepackage{eufrak}
\newcommand{\C}{{\EuFrak C}}

\newcommand{\bi}{\begin{itemize}}
\newcommand{\ei}{\end{itemize}}

\def\Ind#1#2{#1\setbox0=\hbox{$#1x$}\kern\wd0\hbox to 0pt{\hss$#1\mid$\hss}
\lower.9\ht0\hbox to 0pt{\hss$#1\smile$\hss}\kern\wd0}
\def\ind{\mathop{\mathpalette\Ind{}}}
\def\Notind#1#2{#1\setbox0=\hbox{$#1x$}\kern\wd0\hbox to 0pt{\mathchardef
\nn=12854\hss$#1\nn$\kern1.4\wd0\hss}\hbox to
0pt{\hss$#1\mid$\hss}\lower.9\ht0 \hbox to
0pt{\hss$#1\smile$\hss}\kern\wd0}
\def\nind{\mathop{\mathpalette\Notind{}}}

\newtheorem{theorem}{Theorem}[section]
\newtheorem{lemma}[theorem]{Lemma}

\newtheorem{corollary}[theorem]{Corollary}
\newtheorem{proposition}[theorem]{Proposition}
\newtheorem{definition}[theorem]{Definition}
\newtheorem{remark}[theorem]{Remark}
\newtheorem{conjecture}[theorem]{Conjecture}

\newtheorem{problem}[theorem]{Problem}

\title{On stable fields and weight}
\author{Krzysztof Krupi\'nski and Anand Pillay}
\date{September 30th 2009}

\begin{document}
\maketitle
\begin{abstract}
We prove that if $K$ is a (infinite) stable field whose generic type has weight $1$, then $K$ is separably closed. We also obtain some partial results about stable groups and fields whose generic type has finite weight, as well as about strongly stable fields (where by definition {\em all} types have finite weight).
\end{abstract}

\footnotetext{Both authors were supported by the EPSRC grant EP/F009712/1.}

\section{Introduction}
An important aspect  of ``algebraic'' model theory is to uncover the algebraic consequences for structures such as groups, rings, fields,  of abstract model-theoretic properties such as categoricity, stability, simplicity, and so on. Among the first results in the area was Macintyre's theorem \cite{ma} that an infinite field $K$ whose first order theory is $\omega$-stable, is algebraically closed. This was subsequently generalized to superstable fields \cite{CS}, and the proof also works for stable fields with `semiregular generic type''. In all these cases, a suitable rank or dimension is available (for example $U$-rank, or $p$-weight) which one can compute with.

Such methods or tools are on the face of it unavailable in arbitrary stable fields.
Nevertheless a longstanding conjecture is that any infinite field whose first order theory is stable is separably closed.

In the current paper, we discuss this conjecture,  but under additional assumptions on the ``weight'' of the generic type of $K$. We discuss later our motivation.

We refer the reader to \cite{pillay-book} for more details on stability theory, and to \cite{poizat-book} for more details on stable groups and fields. In particular, we assume familiarity with the notion of ``generic type'' of a group, and with the fact that a stable field has a unique generic type.

When we talk about a group $G$ or field $K$ as a first order structure, we mean that the group/field is endowed not only with its algebraic operations, but possible additional relations. When we say for example that $G$ is stable, we mean $Th(G)$ is stable, and by convention $G$ is assumed to be a ``monster model'' or very saturated model of its theory.

The cardinalities of subsets that we are working with are assumed to be smaller than the degree of saturation. And all complete types we mention are finitary, namely are types in finitely many free variables.

\begin{definition}
Let $\C$ be a monster model of a stable theory $T$, $A \subseteq \C$ and $p \in S(A)$. The weight of $p$ is defined as the supremum of the set of cardinalities $\kappa$ for which there exists a non-forking extension $q=tp(a/B) \in S(B)$ of $p$ and a $B$-independent sequence $(b_i: i <\kappa)$ such that $a \nind_A b_i$ for every $i<\kappa$.
\end{definition}

\begin{definition} We say that a stable group $G$ has weight $\alpha$ (symbolically $w(G)=\alpha$) if every (some) generic type of $G$ has weight $\alpha$.
\end{definition}

The weight of a type $p$ is always bounded by $|T|$. In a superstable theory, all types have finite weight; in fact, any type $p$ of weight $n$ will be domination equivalent to a product of $n$ regular types, and any regular type has weight $1$. As mentioned above, there is a machinery ($p$-simplicity, $p$-semiregularity,...) around working close to a regular type in a stable theory, which enables one to prove results as in the superstable case, under the assumption of the existence of enough, or suitable, regular types. However, in a general stable theory, a type can have finite weight or even weight $1$ without being nonorthogonal to a regular type. An important example for the current paper is a separably closed field $K$ of infinite Ershov invariant (or degree of imperfection). It was proved in \cite[Partie IV]{de} that the generic type of $K$ has weight $1$. However, this generic type is not regular (and is, in fact, orthogonal to all regular types). If $K$ is a superstable group, or more generally, a stable group with semiregular generic, then the additivity properties of $U$-rank or $p$-weight can be used to prove what we might call a weak ``exchange property'' for generics: if $g\in G$ is generic (over some fixed set of parameters), $h\in G$ and $g\in acl(h)$ (where $acl(-)$ refers to algebraic closure in the structure $G$ in the model-theoretic sense), then $h$ is also generic. This key property is behind the proofs that for example a superstable field is algebraically closed. But it fails in {\em any} separably closed non perfect field $K$: if $p$ is the characteristic and $g\in K$ is generic, then $g$ is algebraic over $g^{p}$, but $g^{p}$ is not generic. In particular,
it fails in the infinite imperfection degree case but where nevertheless the generic type has weight $1$.

Bearing in mind this example, the strongest conjecture we can make about stable fields of finite weight is:

\begin{conjecture}\label{main conjecture}
Every infinite stable field of finite weight is separably closed.
\end{conjecture}

In this paper, we prove the above conjecture for stable fields of weight 1. We also establish some partial results for stable fields of  finite weight.

Although Conjecture \ref{main conjecture} is interesting in its own right, it is worth giving some motivation.
Shelah recently introduced {\em strongly dependent theories}  as a kind of counterpart of superstable theories, but in the NIP context, and he asked about the structure of strongly dependent fields \cite{sh}. Actually this strong dependence condition turns out to be something like a ``finite weight'' assumption. In fact, assuming stability, strong dependence of $T$ amounts precisely to saying that all types have finite weight \cite[Corollary 9]{ad}. We call strongly dependent stable theories {\em strongly stable}. So, to understand in a meaningful way strongly dependent fields, we would at least have to have some techniques to use a ``finiteness of weight'' hypothesis in the stable case. Thus we were naturally led to ask whether appropriate weight assumptions on the generic type of a stable field could have structural-algebraic consequences.

In section 1, we prove the main theorem (Theorem 1.7) that stable fields of weight $1$ are separably closed. The key lemma shows that we do obtain a kind of weak exchange property for generics under a weight $1$ assumption, but with model-theoretic $acl(-)$ replaced by field-theoretic separable algebraic closure (see Lemma 1.3 and Corollary 1.4).

In section 2, we obtain other partial results around Conjecture \ref{main conjecture}, as well as pointing out in Proposition \ref{perfect} that strongly stable fields are perfect.

The first author is grateful to Frank Wagner for sharing useful ideas.

\section{Stable fields of weight 1}

Let us start from a very important, basic observation. It is essentially contained in Proposition 2.8 of \cite{pillay-freegroup}, but we give a complete proof.

\begin{remark}\label{product of non generics}
Let $G$ be a stable group of weight 1. Then, for an arbitrary set $A$, if $a$ and $b$ are non-generics over $A$, so is the product $a\cdot b$.
\end{remark}
{\em Proof.} Suppose for a contradiction that $a\cdot b$ is generic over $A$. Choose $g$ generic over $A,a,b$. Then $g\cdot a$ is of course generic over $A$. We also have that $g \ind_A a\cdot b$, and so $g \ind_A g\cdot a\cdot b$. On the other hand, $g\cdot a \nind_A g $ (otherwise $g$ is generic over $A,g\cdot a$, so $a=g^{-1}\cdot g\cdot a$ is generic over $A$, a contradiction) and $g\cdot a \nind_A g\cdot a\cdot b$ (otherwise $g\cdot a$ is generic over $A,g\cdot a \cdot b$, so $b=(g\cdot a)^{-1}\cdot g\cdot a\cdot b$ is generic over $A$, a contradiction). Hence, $w(g\cdot a/A) >1$, and so $w(G)>1$, a contradiction.\hfill $\blacksquare$\\

From the above Remark, we get the following
\begin{corollary}

In a stable field $K$ of weight 1, for any $A\subseteq K$, both the sum and the product of two non-generics over $A$ are non-generic over $A$.
\end{corollary}

From now on, in this section, $K$ will be a stable field satisfying the conclusion of the above corollary.
By $p$ we will denote the characteristic of $K$ and by $\mathbb{F}_p$ the prime subfiled of $K$. Also, in the remainder of this section, when we speak of an element of a field being (separably) algebraic over a subfield, we mean of course in the field-theoretic sense.

The following lemma is essential for the proof of the main result.

\begin{lemma}\label{main lemma}
Let $A$ be a subset and $g,h_1,\dots,h_m$ elements of $K$. Suppose $g$ is generic over $A$ and separably algebraic over $\mathbb{F}_p(A,h_1\dots,h_m)$. Then,  $h_i$ is generic over $A$ for some $i\in \{1,\dots,m\}$.
\end{lemma}
{\em Proof.} Put $\overline{h}=(h_1,\dots,h_m)$. Let
$$P(x)=x^n+ R_{n-1}(A,\overline{h})x^{n-1}+\dots+R_0(A,\overline{h})$$
be the minimal polynomial of $g$ over $\mathbb{F}_p(A,\overline{h})$. So, $R_i(A,\overline{y})$'s are rational functions in $\overline{y}$ over $\mathbb{F}_p(A)$, and $P$ is separable.  The proof will be by induction on $n$.

First, consider the base induction step, i.e. $n=1$. We have $g=-R_0(A,\overline{h})$. We can write $R_0(A,\overline{y})=Q(A,\overline{y})/T(A,\overline{y})$, where $Q(A,\overline{y})=\sum a_{i_1,\dots,i_m}y_1^{i_1}\dots y_m^{i_m}$ and $T(A,\overline{y})=\sum b_{j_1,\dots,j_m}y_1^{j_1}\dots y_m^{j_m}$ for some $a_{i_1,\dots, i_m}, b_{j_1,\dots,j_m} \in \mathbb{F}_p(A)$. Since the quotient $Q(A,\overline{h})/T(A,\overline{h})$ is generic over $A$, either $Q(A,\overline{h})$ or $T(A,\overline{h})$ is generic over $A$. Hence, there are $i_1,\dots,i_m$ such that $a_{i_1\dots,i_m}h_1^{i_1}\dots h_{i_m}^{i_m}$ or $b_{i_1\dots,i_m}h_1^{i_1}\dots h_{i_m}^{i_m}$  is generic over $A$. As $a_{i_1,\dots,i_m}$ and $b_{i_1,\dots,i_m}$ are not generic over $A$, we get that one of the $h_i$'s must be generic over $A$, which completes the base induction step.

Now, we turn to the induction step. So, assume that $n>1$ and that the lemma is true for elements whose minimal polynomial has degree smaller than $n$.\\[1mm]
{\bf CASE 1} $p \mid n$.

Since $P(x)$ is separable and irreducible over $\mathbb{F}_p(A)$, there is $1\leq j \leq n-1$ such that $p \nmid j$ and $R_j(A,\overline{h}) \ne 0$.

Take $g_0$ generic over $A,g$. Then, $gg_0$ is generic over $A$. So, $gg_0 \equiv_A g$. Thus, there is $\overline{h'}=(h_1',\dots,h_m') \equiv_A (h_1,\dots,h_m)$ such that
$$(gg_0)^n + R_{n-1}(A,\overline{h'})(gg_0)^{n-1}+\dots +R_0(A,\overline{h'})=0.$$
Put
$$Q(x)=x^n + \frac{R_{n-1}(A,\overline{h'})}{g_0}x^{n-1}+ \dots + \frac{R_{0}(A,\overline{h'})}{g_0^n} \in \mathbb{F}_p(A,g_0,\overline{h'})[x].$$
Let
$$W(x)=Q(x)-P(x) \in \mathbb{F}_p(A,g_0,\overline{h},\overline{h'})[x].$$
We see that $Q(g)=0$, so $W(g)=0$. Moreover,
$$
\begin{array}{ll}
W(x) = & \left( \frac{R_{n-1}(A,\overline{h'})}{g_0} - R_{n-1}(A,\overline{h})\right) x^{n-1} + \dots +\\
&+ \left( \frac{R_{j}(A,\overline{h'})}{g_0^{n-j}} - R_{j}(A,\overline{h})\right) x^{j} + \dots + \left( \frac{R_{0}(A,\overline{h'})}{g_0^n} - R_{0}(A,\overline{h})\right).
\end{array}$$
{\bf Subcase A} $\frac{R_{j}(A,\overline{h'})}{g_0^{n-j}} - R_{j}(A,\overline{h}) =0$.

Then, $g_0^{n-j}-\frac{R_j(A,\overline{h'})}{R_j(A,\overline{h})}=0$. Since $1 \leq n-j<n$ and $p \nmid n-j$, we see that $g_0$ is separably algebraic over $\mathbb{F}_p(A, \overline{h}, \overline{h'})$ and the degree of the minimal polynomial of $g_0$ over $\mathbb{F}_p(A, \overline{h}, \overline{h'})$ is smaller than $n$. Moreover, $g_0$ is generic over $A$. Hence, by the induction hypothesis, there is $i$ such that $h_i$ or $h_i'$ is generic over $A$. But $h_i' \equiv_A h_i$. So, $h_i$ is generic over $A$.\\[1mm]
{\bf Subcase B}  $\frac{R_{j}(A,\overline{h'})}{g_0^{n-j}} - R_{j}(A,\overline{h}) \ne 0$.

We have that $W(g)=0$, $1\leq deg(W) \leq n-1$, and we know that $g$ is separably algebraic over $\mathbb{F}_p(A,g_0,\overline{h},\overline{h'})$. So, the degree of the minimal polynomial of $g$ over this field is smaller than $n$. Moreover, $g$ is generic over $A,g_0$. Hence, by the induction hypothesis,
there is $i$ such that $h_i$ or $h_i'$ is generic over $A,g_0$, so also over $A$.
As $h_i \equiv_A h_i'$, we conclude that $h_i$ is generic over $A$.\\[1mm]
{\bf CASE 2} $p \nmid n$.

Once again, take $g_0$ generic over $A,g$. Then, $g+g_0$ is generic over $A$. So, $g \equiv_A g+g_0$. Thus, there is $\overline{h'} =(h_1'\dots,h_m') \equiv_A \overline{h}$ such that
$$(g+g_0)^n + R_{n-1}(A,\overline{h'})(g+g_0)^{n-1}+\dots +R_0(A,\overline{h'})=0.$$
Put
$$Q(x)=(x+g_0)^n + R_{n-1}(A,\overline{h'})(x+g_0)^{n-1}+ \dots + R_{0}(A,\overline{h'}) \in \mathbb{F}_p(A,g_0,\overline{h'})[x].$$
Let
$$W(x)=Q(x)-P(x) \in \mathbb{F}_p(A,g_0,\overline{h},\overline{h'})[x].$$
We see that $W(g)=0$ and
$$W(x)=(ng_0+R_{n-1}(A,\overline{h'}) - R_{n-1}(A,\overline{h}))x^{n-1}+W_1(x),$$
where $W_1(x) \in \mathbb{F}_p(A,g_0,\overline{h},\overline{h'})[x]$ is of degree smaller than $n-1$.\\[1mm]
{\bf Subcase A} $ng_0+R_{n-1}(A,\overline{h'}) - R_{n-1}(A,\overline{h})=0$.

Since $p \nmid n$ and $g_0$ is generic over $A$, by the base induction step, we get that there is $i$ such that $h_i$ or $h_i'$ is generic over $A$. So, $h_i$ is generic over $A$.\\[1mm]
{\bf Subcase B} $ng_0+R_{n-1}(A,\overline{h'}) - R_{n-1}(A,\overline{h})\ne 0$.

Then, $W(g)=0$, $deg(W)=n-1\geq 1$, and we know that $g$ is separably algebraic over $\mathbb{F}_p(A,g_0,\overline{h},\overline{h'})$. So, the degree of the minimal polynomial of $g$ over this field is smaller than $n$.  Moreover, $g$ is generic over $A,g_0$. Hence, we finish using the induction hypothesis as in Subcase B of Case 1. \hfill $\blacksquare$\\

Notice that if the characteristic of $K$ equals 0, then Case 1 does not hold, and so it is enough to apply the argument from Case 2 to prove Lemma \ref{main lemma}.

Let us formulate Lemma \ref{main lemma} in the case $m=1$ as a corollary. 

\begin{corollary}\label{R-condition}
Let $A$ be a subset and $g,h$ elements of $K$. Suppose $g$ is generic over $A$ and separably algebraic over $\mathbb{F}_p(A,h)$. Then, $h$ is generic over $A$.
\end{corollary}

\begin{lemma}\label{R for tuples}
Let $A$ be a subset of $K$ and $g_1,\dots, g_m$ independent generics over $A$. Suppose $h_1,\dots,h_m$ are such that the elements $g_1,\dots,g_m$ are separably algebraic over $\mathbb{F}_p(A,h_1,\dots,h_m)$. Then, $h_1,\dots,h_m$ are independent generics over $A$.
\end{lemma}
{\em Proof.} The proof is by induction on $m$. For $m=1$, the conclusion follows from Corollary \ref{R-condition}.

Let us do the induction step. By the assumption, $g_m$ is generic over $A,g_{<m}$ and it is separably algebraic over $\mathbb{F}_p(A,g_{<m},\overline{h})$. Hence, by Lemma \ref{main lemma}, there is $i$ such that $h_{i}$ is generic over $A,g_{<m}$. Therefore, $h_i,g_1,\dots,g_{m-1}$ are independent generics over $A$.

Put $A'=A \cup \{h_i\}$. We see that $g_1,\dots,g_{m-1}$ are independent generics over $A'$ and they are separably algebraic over $\mathbb{F}_p(A',h_{\ne i})$. So, by the induction hypothesis, $h_1,\dots,h_{i-1},h_{i+1},\dots,h_m$ are also independent generics over $A'$. We finish using the fact that $h_i$ is generic over $A$. \hfill $\blacksquare$

\begin{corollary}\label{symmetric functions}
Let $A$ be a subset of $K$ and $a_0,\dots,a_{n-1}$ independent generics over $A$. Then:\\
(i) the elementary symmetric functions in $a_0,\dots,a_{n-1}$ are independent generics over $A$,\\
(ii) the polynomial $x^n+a_{n-1}x^{n-1}+\dots + a_0$ has $n$ distinct roots in $K$.
\end{corollary}
{\em Proof.}
(i) Let $s_0,\dots,s_{n-1}$ be the elementary symmetric functions in $\overline{a}$. We have that $a_0,\dots,a_{n-1}$ are pairwise distinct solutions to $x^n-s_{n-1}x^{n-1}+\dots + (-1)^ns_0$. Hence, $a_0,\dots,a_{n-1}$ are separably algebraic over $\mathbb{F}_p(A,s_0,\dots,s_{n-1})$. So, by Lemma \ref{R for tuples}, we get that $s_0,\dots,s_{n-1}$ are independent generics over $A$.\\
(ii) It follows from (i) and the uniqueness of the generic type. \hfill $\blacksquare$\\

With the above lemmas and corollaries, we can now prove our main result, by adapting the proof of \cite[Proposition 5.2]{pi}.

\begin{theorem}\label{main theorem}
Each stable field of weight 1 is separably closed.
\end{theorem}
{\em Proof.} As usual, $K$ is our stable field of weight 1 and $p$ is its characteristic. Suppose for a contradiction that there is $\alpha \in K^{sep} \setminus K$. Let $P(x)=x^n+a_{n-1}x^{n-1}+\dots + a_0$ be the minimal polynomial of $\alpha$ over $K$. Since $\alpha \in K^{sep}$, $P(x)$ has $n$ different roots $\alpha_1,\dots, \alpha_n$ in $K^{sep}$.
Choose
\begin{enumerate}
\item[(i)] $t_0,\dots, t_{n-1}$ independent generics over $a_0,\dots, a_{n-1}$.
\end{enumerate}
Define
$$r_i=t_0+t_1\alpha_i+\dots + t_{n-1}\alpha_i^{n-1}$$
for $i=1,\dots,n$. Let $s_0, \dots,s_{n-1}$ be the elementary symmetric function in $r_1,\dots,r_n$. Then, $s_0,\dots,s_{n-1} \in K$ because they are fixed by every element of $Gal(K^{sep}/K)$. We claim that
\begin{enumerate}
\item[(ii)] $r_1,\dots,r_n$ are separably algebraic over $\mathbb{F}_p(s_0,\dots,s_{n-1})$.
\end{enumerate}
We have that $r_1,\dots, r_n$ are the roots of $x^n-s_{n-1}x^{n-1}+\dots +(-1)^ns_0$. So, in order to prove (ii), it is enough to show that $r_i \ne r_j$ whenever $i\ne j$. Suppose for a contradiction that there are $i \ne j$ such that $r_i=r_j$. Then, $$t_1(\alpha_i-\alpha_j)+\dots +t_{n-1}(\alpha_i^{n-1} - \alpha_j^{n-1})=0.$$ So, $t_1$ is algebraic over $\mathbb{F}_p(\alpha_i,\alpha_j,t_2,\dots,t_{n-1})$ and so over $\mathbb{F}_p(a_0,\dots,a_{n-1},t_2,\dots,t_{n-1})$, which contradicts (i).

Since the matrix
$$\left( \begin{array}{llll}
1 & \alpha_1 & \dots & \alpha_1^{n-1}\\
\vdots & \vdots & \vdots & \vdots\\
1 & \alpha_n & \dots & \alpha_n^{n-1}
\end{array}\right)$$
is invertible, we see that $t_0,\dots,t_{n-1} \in \mathbb{F}_p(\alpha_1,\dots,\alpha_n,r_1,\dots,r_n)$. On the other hand, $\alpha_1,\dots,\alpha_{n-1} \in \mathbb{F}_p(a_0,\dots,a_{n-1})^{sep}$. Thus, by (ii),
\begin{enumerate}
\item[(iii)] $t_0,\dots,t_{n-1}$ are separably algebraic over $\mathbb{F}_p(a_0,\dots,a_{n-1},s_0,\dots,s_{n-1})$.
\end{enumerate}
By (i), (iii) and Lemma 1.5, we see that $s_0,\dots,s_{n-1}$ are independent generics over $a_0,\dots,a_{n-1}$. So, in virtue of Corollary \ref{symmetric functions}(ii), all $r_i$'s belong to $K$. Thus, the degree of the minimal polynomial of $\alpha$  over $K$ is smaller than $n$, a contradiction. \hfill $\blacksquare$\\

Recall that by \cite{de}, we know that the weight of a separably closed field of infinite Ershov invariant is 1. Thus, Theorem \ref{main theorem} is in a sense best possible.

\section{Stable fields of finite weight}

\begin{proposition}\label{groups of finite weight} Let $G$ be any stable commutative group (written multiplicatively) of finite weight. Then for all but finitely many primes $q$, $G^q$ has finite index in $G$. More precisely, if $w(G)=w<\omega$, then there are at most $w$ many primes $q$ such that $[G:G^q]$ is infinite.
\end{proposition}
{\em Proof.}
Choose an independent sequence $(a_n)_{n \in \omega}$ of generics in $G$. Assume $w(G)=w < \omega$. Suppose for a contradiction that there are $w+1$ primes $p_1,\dots,p_{w+1}$ such that $[G:G^{p_i}]$ are infinite. It follows that $G^{p_i}$ are not generic.

Define a sequence $(k_1,\dots,k_{w+1})$ of natural numbers by

$$\left\{
\begin{array}{lrl}
k_1 & = & p_1+1,\\
 k_{i} & = & (p_1\dots p_{i-1})^{p_i -1}\; \, \mbox{for} \; 2\leq i \leq w+1.
\end{array} \right.
$$
Then, $p_i|k_i-1$ for any $i=1,\dots,w+1$, and $p_i|k_j$ for any $i<j$.

Put $g=a_0a_1^{k_1}\dots a_{w+1}^{k_{w+1}}$, and define a sequence $(g_i)_{1\leq i \leq w+1}$ of elements of $G$ by

$$\left\{
\begin{array}{lrl}
g_1 & = & a_0a_1,\\
g_i & = & a_0a_1^{k_1}\dots a_{i-1}^{k_{i-1}}a_i \; \, \mbox{for} \; 2 \leq i \leq w+1.
\end{array} \right.
$$\\[3mm]
{\bf Claim}
(i) $g$ is generic.\\
(ii) $g_i \ind g_j$ for any $i \ne j$.\\
(iii) $g \nind g_{i}$ for every $i=1,\dots,w+1$.\\[3mm]
{\em Proof of Claim.}
(i) It follows from the fact that $a_0$ is generic over $a_{>0}$.\\
(ii) Assume $j >i$. Since $a_j$ is generic over $a_{<j}$,
$$g_j=a_0a_1^{k_1}\dots a_{j-1}^{k_{j-1}}a_j \ind a_0a_1^{k_1}\dots a_{i-1}^{k_{i-1}}a_i=g_i.$$
(iii) We have $g_i^{-1}g=a_i^{k_i -1}a_{i+1}^{k_{i+1}}\dots a_{w+1}^{k_{w+1}}$. Since $p_i|k_i-1$ and $p_i|k_j$ for any $j>i$, we get that $g_i^{-1}g \in G^{p_i}$.
Suppose for a contradiction that $g \ind g_{i}$. Then, by (i), $g_{i}^{-1}g$ is generic, and hence $G^{p_i}$ is generic, a contradiction. \hfill $\square$\\

By the Claim, $w(K)\geq w+1$, a contradiction. \hfill $\blacksquare$

\begin{corollary}\label{almost all} Let $K$ be any infinite stable field of finite weight. Then for all but finitely many primes $q$, $K^q=K$. More precisely, if $w(K)=w<\omega$, then there are at most $w$ many primes $q$ for which $K^q \ne K$.
\end{corollary}
{\em Proof.} It an immediate consequence of Proposition \ref{groups of finite weight} and the fact that stable fields are multiplicatively connected.\hfill $\blacksquare$\\

As the weight of a separably closed field of infinite Ershov invariant is 1, we cannot expect to strengthen Corollary \ref{almost all} to get that for every prime $q$, $K^q=K$. However, one can hope to prove that for every prime $q$ different from the characteristic, $K^q=K$. In fact, this would imply Conjecture \ref{main conjecture}. To see this, one should apply  Macintyre's proof \cite{ma} using the fact that a finite extension of a stable field of finite weight remains stable of finite weight and the fact that stable fields are closed under Artin-Schreier extensions \cite{KSW}.


As was mentioned in the introduction, a separably closed field of infinite Ershov invariant is an example of stable field of finite weight which is not strongly stable, i.e. there is a finitary type in it of infinite weight. This follows from the next proposition.

\begin{proposition}\label{perfect}
An infinite strongly stable field is perfect.
\end{proposition}
{\em Proof.}  Let $p$ be the characteristic of $K$. Assume $p>0$, and suppose for a contradiction that $K^p \ne K$. Then, there are $b_1,b_2 \in K$ linearly independent over $K^p$. Choose a Morley sequence $(a_i)_{i \in \omega}$ in the generic type over $b_1,b_2$.

By compactness, one can find $a\in K$ for which there is a sequence $(c_i)_{i \in \omega}$ of elements of $K$ such that $c_0=a$ and for every $i$, $c_i=b_1c_{i+1}^p+b_2a_i^p$.

Since $b_1,b_2$ are linearly independent over $K^p$, we get that $a_i \in dcl(b_1,b_2,a)$ for every $i$. So $a \nind_{b_1,b_2} a_i$ for every $i$. On the other hand, $(a_i)_{i \in \omega}$ was chosen to be independent over $b_1,b_2$. So $w(a/b_1,b_2)$ is infinite, and hence $K$ is not strongly stable, a contradiction. \hfill $\blacksquare$\\

The above proposition together with Theorem \ref{main theorem} yield the following corollary.

\begin{corollary}
Strongly stable fields of weight 1 are algebraically closed.
\end{corollary}

The next observation says that if we assume that the degree of imperfection is finite, then the conclusion of Proposition \ref{perfect} holds under the weaker assumption of being of finite weight.

\begin{proposition}
An infinite stable field of finite weight and of finite degree of imperfection is perfect.
\end{proposition}
{\em Proof.} Let $p>0$ be the characteristic of $K$. Suppose for a contradiction that $K^p \ne K$. Since the degree of imperfection is finite, there is a finite basis $\{ b_1,\dots,b_n\}$ of $K$ over $K^p$.

The map $f: K \to K^{\times n}$ given by $f(a)=(f_1(a),\dots,f_n(a))$ where $a=b_1f_1(a)^p +\dots + b_nf_n(a)^p$ is a group automorphism definable over $\{b_1,\dots,b_n\}$. So for any $A\subseteq K$, if $a$ is generic over $A,b_1,\dots,b_n$, then $f(a)$ is a sequence of independent generics over $A,b_1,\dots,b_n$. For $\eta \in \{1,\dots,n\}^l$ and $x \in K$, we put $f_\eta(x)=(f_{\eta(l-1)}\circ \dots \circ f_{\eta(0)})(x)$.

Let $a$ be generic over $b_1,\dots,b_n$. By an easy induction, we get that $(f_{1^ki}(a): k \geq 0, 1<i\leq n)$ is an infinite collection of independent generics over $b_1,\dots,b_n$ ($1^ki$ denotes the sequence consisting of $k$ many 1's followed by $i$). Moreover, every $f_{1^ki}(a)$ belongs to $dcl(a,b_1,\dots,b_n)$. So, $a \nind_{b_1,\dots,b_n} f_{1^ki}(a)$. We conclude that $w(a/b_1,\dots,b_n)$ is infinite. Thus, $w(K)$ is infinite, a contradiction. \hfill $\blacksquare$\\

We complete the paper with a couple of questions and conjectures related to the notions and techniques introduced here. We did not give much thought to the first one, but we are rather curious and there could be a simple construction.

\begin{problem} Construct an algebraically closed field $K$ with additional structure such that $Th(K)$ is stable and the generic type of $K$ has weight $1$ but is not regular.
\end{problem}

\begin{conjecture} Let $K$ be a field with additional structure which is stable (and saturated).  Then, the following are equivalent:
\newline
(1) $K$ is separably closed,
\newline
(2)  For any small subfield $k<K$, $n < \omega$, and $a_{1},\dots,a_{n},b_{1},\dots,b_{n}\in K$, IF $a_{1},\dots,a_{n}$ are independent generics over $k$, and each $a_{i}$ is separably algebraic over $k(b_{1},\dots,b_{n})$ (of course, in the field-theoretic sense), THEN $b_{1},\dots,b_{n}$ are independent generics over $k$.
\end{conjecture}

Note that (2) is precisely the statement of Lemma 1.5, and the proof of Theorem 1.7 shows that (2) implies (1).

Here is a version for ``algebraically closed'' rather than ``separably closed''.
\begin{conjecture} Let $K$ be a field with additional structure which is stable (and saturated). Then, the following are equivalent:
\newline
(1') $K$ is algebraically closed,
\newline
(2') For any small subfield $k<K$ and $a,b\in K$, IF $a$ is generic over $k$, and $a$ is algebraic over $k(b)$ in the field-theoretic sense, THEN $b$ is generic over $k$.
\end{conjecture}

In fact, it is not hard to show that (2') implies (2''): if $a_{1},\dots,a_{n}$ are generic independent over $k$ and contained in $k(b_{1},\dots,b_{n})^{alg}$, then $b_{1},\dots,b_{n}$ are generic independent over $k$. And by the standard argument (as in the proof of 1.7), one deduces from (2'') that $K$ is algebraically closed. The converse (1') implies (2') looks attractive, and concerns some kind of uniqueness of ``generic types'' on irreducible plane curves in stable expansions of algebraically closed fields.

In any case, the point of Conjectures 2.7 and 2.8 is that the kind of methods in the proof of Theorem 1.7 are not only sufficient but should also be necessary.

\noindent
Krzysztof Krupi\'nski\\
Instytut Matematyczny Uniwersytetu Wroc\l awskiego\\
pl. Grunwaldzki 2/4, 50-384 Wroc\l aw, Poland.\\
e-mail: kkrup@math.uni.wroc.pl\\[5mm]
Anand Pillay\\
Department of Pure Mathematics\\
University of Leeds\\
Leeds, LS2 9JT.\\
e-mail: pillay@maths.leed.ac.uk

\end{document}